\documentclass[11pt]{article}
\usepackage{amssymb,latexsym,amsmath,amsbsy,amsthm,amsxtra,amsgen}
\newfont{\msbm}{msbm10 at 11pt}

\def\mn{\bigskip\noindent}
\def\ZZ{{\bf Z}}
\def\la{\lambda}
\def\ee{\epsilon}

\newtheorem{Theo}{Theorem}
\newtheorem{Lemma}[Theo]{Lemma}

\newtheorem{Prop}[Theo]{Proposition}

\begin{document}
\title{Spatial and non-spatial stochastic models for immune response}
\author{by Rinaldo B. Schinazi\thanks{Supported in part by a NSA grant} \:
and Jason Schweinsberg\thanks{Supported in part by NSF Grant DMS-0504882} \\ \\
University of Colorado at Colorado Springs and University of California at San Diego}
\maketitle

\begin{abstract}
We propose both spatial and non-spatial stochastic models for pathogen dynamics in the presence of an immune response. One of our spatial models shows that, at least in theory, a pathogen may escape the immune system thanks to its high mutation probability alone.  While one of our non-spatial models also exhibits this behavior, another behaves quite differently from the corresponding spatial model.
\end{abstract}


\section{Introduction and results}

We study some simple mathematical models designed to test the following hypothesis: can a pathogen escape the immune system only because of its high probability of mutation?  We propose both spatial and non-spatial models.  In all of our models, we assume that pathogens can mutate, leading to the appearance of new types of pathogens.  We also assume that the immune system is able to get rid of all the pathogens of a given type at once but that it recognizes only one type at a time. 

\subsection{Non-spatial models}

For our non-spatial models, we start with a single pathogen at time zero.  Each pathogen gives birth to a new pathogen at rate $\lambda$.  When a new pathogen is born, it has the same type as its parent with probability $1 - r$.  With probability $r$, a mutation occurs, and the new pathogen has a different type from all previously observed pathogens.  For convenience, we say that the pathogen present at time zero has type $1$, and the $k$th type to appear will be called type $k$.  Note that we assume the birth rate $\la$ to be the same for all types and we therefore ignore selection pressures.

If there are no deaths, then this is a model of a Yule process with infinitely many types, which goes back to Yule (1925).  The model with no deaths was studied recently by Durrett and Schweinsberg (2005), who focused on the joint distribution of the number of pathogens of each type.  Here we assume that the response of the immune system can eliminate pathogens of a given type.  We propose three different models for the behavior of the immune system. In all three models the pathogens give birth and mutate as described above. Each model corresponds to a different immune response.

\mn {\bf Model 1}:  At times of a rate $1$ Poisson process, a death event occurs.  When there is a death event, if there are $k$ types of pathogens alive, then one of the types is chosen at random, each with probability $1/k$, and all pathogens of that type are simultaneously killed.

\mn {\bf Model 2}:  When a new type appears in the population, it survives for an exponential amount of time with mean $1$, independently of all the other types.  All pathogens of the type are killed simultaneously.

\mn {\bf Model 3}:  Each pathogen is born with a mean 1 exponential clock. When the clock goes off the pathogen is killed as well as all the pathogens of the same type. 

\bigskip
To understand the models better, note that if there are $k$ types and $N$ total pathogens, then the total rate of death events is $1$ in Model 1, $k$ in Model 2, and $N$ in Model 3.  Also, if there are $n_i$ pathogens of type $i$, the rate at which type $i$ is being killed is $1/k$ in Model 1, $1$ in Model 2, and $n_i$ in Model 3.  Thus, in models 1 and 2, the rate at which a type is killed does not depend on the number of pathogens with that type, but in Model 3, types that have large numbers of pathogens are more likely to be targeted by the immune system and eliminated.  Model 2 is similar to random graph models studied by Chung and Lu (2004) and Cooper, Frieze, and Vera (2004) with preferential attachment (corresponding to births) and vertex deletion (corresponding to deaths).

With all three models, there is a positive probability that the immune system will succeed in eliminating the pathogens, as the first pathogen could die before it has any offspring.  Our main result for the non-spatial models is the following theorem, which specifies the values of $r$ and $\lambda$ for which there is a positive probability that the pathogens survive, meaning that for all $t > 0$, there is at least one pathogen alive at time $t$.

\begin{Theo} 
Assume $\lambda > 0$ and $r > 0$.
\begin{enumerate}
\item For Model 1, the pathogens survive with positive probability.

\item For Model 2, the pathogens survive with positive probability if and only if $\lambda > 1$.

\item For Model 3, the pathogens survive with positive probability if and only if $r > 1/\lambda$.
\end{enumerate}
\label{nonspatialthm}
\end{Theo}

Thus, the three models produce very different behavior.  In Model 1, which is the model in which the death rates are lowest, any positive probability of mutation is enough to allow the pathogen to escape the immune system.  For Model 2, whether or not the pathogen can survive depends only on the reproduction rate and not on the mutation rate.  For Model 3, there is a phase transition, in that for fixed $\lambda$, the pathogens can escape the immune system if $r > 1/\lambda$ but not if $r \leq 1/\lambda$.  The proof of Theorem \ref{nonspatialthm} is given in Section 2.

\subsection{Spatial models}

We now introduce three spatial models that correspond to the three non-spatial models above.
Each spatial stochastic model is on the lattice $\ZZ^d$, where the dimension $d$ can be any positive integer.   
Every site of $\ZZ^d$ is either occupied by a pathogen or empty. Each model is started  with a single pathogen at the origin of $\ZZ^d$ and with all other sites empty.

Rules for births and mutations are the same for the three spatial models.
Let $x$ be a site occupied by a pathogen and $y$ be one of its $2d$ nearest neighbors. After a random exponential time with rate $\lambda$, the pathogen on $x$ gives birth on $y$, provided $y$ is empty (if $y$ is occupied nothing happens). With probability $1-r$ the new pathogen on $y$ is of the same type as the parent pathogen on $x$. With probability $r$ the new pathogen is of a different type. We assume that every new type that appears is different from all types that have ever appeared.

In addition to the birth and mutation rates decribed above, the three spatial models, which we call S1, S2 and S3, have the same rules for the immune responses as non-spatial models 1, 2 and 3, respectively. 
We start with a result for Model S1.  The result shows that this model produces the same behavior as the corresponding non-spatial model.

\begin{Theo}
Consider Model S1 on $\ZZ^d$ for $d\geq 1$. For every $\lambda>0$ and $r>0$, the pathogens have a positive probability of surviving.
\label{spatial1thm}
\end{Theo}

We now turn to models S2 and S3. If $r=1$, then every birth gives rise to a new type in models S2 and S3. Since all pathogens are of different types there is only one death at a time in both models. If we ignore the types, the process of occupied sites is the well-known contact process. The contact process has a critical value $\lambda_c$ which depends on the dimension $d$ of the lattice.  If $\lambda \leq \lambda_c$ the pathogens die out, while if $\lambda>\lambda_c$ there is a positive probability that pathogens will survive forever. For more on the contact process, see Liggett (1999). 


\begin{Theo}

\begin{enumerate}
\item Consider Model S2 with $\lambda\leq 1/2d$.  For all $r$ in $[0,1]$ the pathogens die out with probability $1$.

\item Consider Model S3 with $\lambda \leq \lambda_c$.  For all $r$ in [0,1], the pathogens die out with probability $1$.
\end{enumerate}
\label{newthm}
\end{Theo}

\begin{Theo}
Consider Models S2 and S3 on $\ZZ^d$ for $d\geq 1$ with parameters $\lambda$ and $r$.
\begin{enumerate} 
\item For any $\lambda>\lambda_c$, there is an $r_1$ in (0,1) such that if $r<r_1$, then the pathogens die out. 

\item For any $\lambda>\lambda_c$, there is an $r_2$ in (0,1) such that if $r>r_2$, then the pathogens survive with positive probability.
\end{enumerate}
\label{spatial2thm}
\end{Theo}

We conjecture that for both Model S2 and S3 there is a critical value $r_c$ such that the pathogens die out if $r<r_c$ and survive if $r>r_c$. This would follow from our results if we could prove, for instance, that the probability of pathogen survival is increasing in $r$.  However, it is not clear that this is a true statement.

While Models 3 and S3 behave in a similar way, models 2 and S2 are strikingly different. In particular, Model S2 exhibits a phase transition in $r$ (survival of pathogens for large $r$, death for small $r$) and Model 2 does not.

Theorems \ref{newthm} and \ref{spatial2thm} will be proved in Section 3, and Theorem \ref{spatial1thm} will be proved in Section 4.
\subsection{Discussion}
In this paper we propose a new class of models for immune response. Our main assumption is that 
the immune system is able to get rid of all the pathogens of a given type at once. As far as we know 
these are the first models that attempt to mimic the ``central command" nature of the immune system with global killing rules. There is increasing evidence that the immune system has a coordinated global behavior. See in particular Silvestri and Feinberg (2003) who argue that the pathogenesis of AIDS is caused by chronic immune activation rather than direct attack of the HIV virus. In contrast, most of the existing models are predator-prey differential equations models with local killing rules, see De Boer and Perelson (1998), Iwasa et al. (2004) and Nowak and May (2000). There exist also a number of spatial models, in particular cellular automata, to model HIV infections, see Perelson and Weisbuch (1997), Bernaschi and Castiglione (2002) or Zorzenon (1999). However, these models, unlike ours, are quite complex, use local killing rules and are not analyzed rigorously.

Our original motivation for the introduction of our models is to to test the following hypothesis: can the 
immune system be overwhelmed by a particular virus only because of its high probability of mutation?
Ordinary differential equation models have been used to test this hypothesis. In particular, Nowak and May (2000) (see Sections 12.1 and 12.2) introduce models of increasing complexity to get a behavior similar to the behavior exhibited by our simple models 3, S2 and S3 (pathogens die out for small mutation rate and survive for large mutation rate).  
Sasaki (1994) uses a partial differential equation model in which all types of pathogens have the same reproduction rate.  However, his analysis yields results strikingly different from ours. In particular, he finds that the pathogens may survive only if the mutation rate is intermediate. If the mutation rate is too low or too high, the pathogens die out in his model; see in particular his results for the infinite allele model. 

The rest of the paper is devoted to proofs. We are able to give short proofs for all our results except for Theorem 2 when $d=1$. While this is probably not the most biologically significant model we feel that our mathematical analysis is worthwhile. 
The proof that an interacting particle system survives is almost always done by coupling the system to a much simpler one. This is the case in this paper for all spatial models except for Model S2 in $d=1$ for
which we could not find such a coupling. Instead, to prove that the pathogens survive we do a ``pathwise" analysis of our model. This is a rather delicate analysis, but it yields a lot more information about the behavior of the system than the coupling technique does.

\section{Analysis of the non-spatial models}

Theorem \ref{nonspatialthm} can be proved using standard branching process techniques.  We begin with the analysis of Model 1, which can be carried out using a comparision with a birth-death process.

\begin{proof}[Proof of part 1 of Theorem \ref{nonspatialthm}]
Let $X(t)$ be the number of different types of pathogens alive at time $t$.  Thus, $X(0) = 1$, and we need to show that for all $r > 0$, we have $P(X(t) > 0 \mbox{ for all }t) > 0$.  Since $r > 0$, we can choose $N$ such that $N \lambda r > 1$.  Let $T_0 = \inf\{t: X(t) = N\}$.  With positive probability, the pathogen present at time zero gives birth to pathogens of $N$ new types before the first death event.  Therefore, $P(T_0 < \infty) > 0$.  For positive integers $n$, inductively define, on the event that $T_{n-1} < \infty$, the stopping time $T_n = \inf\{t > T_{n-1}: X(t) \neq X(T_{n-1})\}$.  That is, $T_n$ is the first time after $T_{n-1}$ that either a pathogen of a new mutant type is born, in which case $X(T_n) = X(T_{n-1}) + 1$, or one of the types is eliminated, in which case $X(T_n) = X(T_{n-1}) - 1$.

Since there must be at least one pathogen of each type, and each pathogen gives birth to pathogens of new types at rate $\lambda r$, at time $t$ the rate at which new types are being born is at least $X(t) \lambda r$.  Therefore, whenever $X(t) \geq N$, the rate at which new types are being born is at least $N \lambda r$.   The rate of death events, which cause a type to be eliminated, is always $1$.  Let $p = N \lambda r/(1 + N \lambda r) > 1/2$.  Then for all $k \geq N$, we have $P(X(T_n) = k+1|X(T_{n-1}) = k) \geq p$.  
Now, consider a birth-death process $(Y_n)_{n=0}^{\infty}$ such that $Y_0 = N$ and, for all integers $k$, we have $P(Y_n = k+1|Y_{n-1} = k) = p$ and $P(Y_n = k-1|Y_{n-1} = k) = 1-p$.  It is well-known, see for 
instance Hoel, Port and Stone (1972), p. 32, that $P(Y_n \geq N \mbox{ for all }n) = (2p-1)/p > 0$.  
Therefore, by comparing the processes $(X(T_n))_{n=0}^{\infty}$ and $(Y_n)_{n=0}^{\infty}$, we see that $P(X(T_n) \geq N \mbox{ for all }n|T_0 < \infty) \geq (2p-1)/p$.  On the event that $T_0 < \infty$ and $X(T_n) \geq N$ for all $n$, we have $X(t) > 0$ for all $t$.  Therefore, the probability that the pathogens survive is at least $P(T_0 < \infty) (2p-1)/p > 0$.
\end{proof}

To study Models 2 and 3, we will construct a tree which keeps track of the genealogy of the different types of pathogens.  Each vertex in the tree will be labeled by a positive integer.  There will be a vertex labeled $k$ if and only if a pathogen of type $k$ is born at some time.  We draw a directed edge from $j$ to $k$ if the first pathogen of type $k$ to be born had a pathogen of type $j$ as its parent.  This construction gives a tree whose root is labeled $1$ because all types of pathogens are descended from the pathogen of type $1$ that is present at time zero.  Since every type is eliminated eventually, we have $X(t) > 0$ for all $t$ if and only if infinitely many different types of pathogens eventually appear or, in other words, if and only if the tree described above has infinitely many vertices.

For Model 1, the rate at which a given type is killed depends on the number of other types present.  However, for Models 2 and 3, the rate at which a type is killed is either constant in the case of Model 2, or depends only on the number of pathogens of the type in the case of Model 3.  Consequently, once the first pathogen of type $k$ is born, the number of mutant offspring born to type $k$ pathogens is independent of how the other types evolve.  Therefore, the tree constructed above is a Galton-Watson tree, and the process survives with positive probability if and only if the mean of the offspring distribution is greater than one, see for instance I.9 in Schinazi (1999).  This observation can be used to prove parts 2 and 3 of Theorem \ref{nonspatialthm}.  We begin with part 3, which is simpler.

\begin{proof}[Proof of part 3 of Theorem \ref{nonspatialthm}]
Whenever there are $n$ pathogens of a given type, the event in which the type is destroyed is happening at rate $n$, while events in which pathogens of the type give birth to offspring of new mutant types are happening at rate $n r \lambda$.  Therefore, regardless of the number of pathogens of the given type, the probability that the type is destroyed before the next birth to a mutant type is $n/(n + n r \lambda) = 1/(1 + r \lambda)$, and the probability that an individual gives birth to a mutant offspring before the type is destroyed is $r \lambda/(1 + r \lambda)$.  Suppose $X$ is the number of types that are offspring of a given type.  Then for $k \geq 0$, $$P(X = k) = \bigg( \frac{r \lambda}{1 + r \lambda} \bigg)^k \bigg( \frac{1}{1 + r \lambda} \bigg) = \frac{(r \lambda)^k}{(1 + r \lambda)^{k+1}}.$$  That is, $X + 1$ has the geometric distribution with parameter $1/(1 + r \lambda)$.  It follows that the mean of the offspring distribution is greater than one if and only if $r > 1/\lambda$.  As discussed above, this is the condition for the process to survive with positive probability.
\end{proof}

\begin{proof}[Proof of part 2 of Theorem \ref{nonspatialthm}]
If $r = 1$, then there is only one individual of each type, so we have births at rate $\lambda$ and deaths of a single individual at rate $1$.  In this case, the result is a standard fact about branching processes.  Now suppose $r < 1$.  Let $m'$ be the number of type $1$ pathogens which are offspring of the initial pathogen.  Note that the type $1$ pathogens evolve like a Yule process with births at rate $\lambda(1 - r)$ until the type $1$ pathogens all die at time $T$, which has an exponential distribution with mean one.  If $Y(t)$ denotes the number of type $1$ pathogens at time $t$, then conditioning on the value of $T$ gives $$m' + 1 = E[Y(T)] = \int_0^{\infty} e^{-t} E[Y(t)] \: dt = \int_0^{\infty} e^{-t} e^{\lambda(1-r)t} \: dt = \int_0^{\infty} e^{-(1 - \lambda(1 - r))t} \: dt.$$
It follows that $m' = \infty$ if $\lambda(1 - r) \geq 1$ and $m' = \lambda(1-r)/(1 - \lambda(1-r))$ if $\lambda(1-r) < 1$.  

Now, let $m$ be the mean number of different types that are offspring of type 1 pathogens.  Because each type $1$ pathogen gives birth to new types at rate $r \lambda$ and to other type $1$ pathogens at rate $(1-r) \lambda$, we must have $m = rm'/(1-r)$.  Therefore, $m = \infty$ if $\lambda(1 - r) \geq 1$ and 
$m = r \lambda/(1 - \lambda(1 - r))$ if $\lambda(1 - r) < 1$.  It follows that $m > 1$ if and only if $\lambda > 1$.  Thus, the process survives with positive probability if and only if $\lambda > 1$.
\end{proof}

\section{Analysis of the second and third spatial models}

In this section, we prove Theorems \ref{newthm} and \ref{spatial2thm}, which pertain to the spatial models S2 and S3.


\begin{proof}[Proof of part 1 of Theorem \ref{newthm}]

Suppose $\lambda \leq 1/2d$.  We may couple Model S2, in which an individual gives birth on each of the $2d$ neighboring sites at rate $\lambda$, with Model 2, in which each individual gives birth at rate $\lambda' = 2d \lambda$.  In both models, each type dies at rate 1, and each pathogen gives birth at rate $2d \lambda$.  However, a birth that occurs in Model 2 will be suppressed in Model S2 if the site on which a pathogen is to give birth is already occupied.  Hence, at any given time, Model 2 has at least as many types as Model S2, and each type in Model 2 has at least as many individuals as the corresponding type in Model S2.

According to Theorem 1.2, if $\lambda' \leq 1$ then the pathogens in Model 2 die out with probability 1. Using the coupling above, one sees that the same is true for Model S2 if $\lambda \leq 1/2d$. 
\end{proof}
\begin{proof}[Proof of part 2 of Theorem \ref{newthm}]

As noted before, when $r=1$ the process of occupied sites is a contact process for models S2 and S3. Since $\lambda \leq \lambda_c$ the pathogens die out. We may couple, site by site, Model S3 with $r<1$ to Model S3 with $r=1$. Deaths occur simultaneously in both models at the same rate 1, but every time there is a death in the process with $r=1$, a single pathogen dies while for the model with $r<1$ all pathogens of the same type die. Birth rates of pathogens are the same. It is easy to see that, with this coupling, the model with $r<1$ has fewer occupied sites than the model with $r=1$, at all times. Since the pathogens die out for Model S3 with $r=1$ they also must die out for S3 with $r<1$.   
\end{proof}

\begin{proof}[Proof of part 1 of Theorem \ref{spatial2thm}]
We write the proof in dimension $d=2$. The same ideas work in any $d\geq 1$.
We start by defining two space-time regions.
$${\cal A}=[-2L,2L]^2\times[0,2T]\qquad {\cal B}=[-L,L]^2\times[T,2T].$$
Note that ${\cal B}$ is nested in ${\cal A}$. 
Let ${\cal C}$ be the boundary of ${\cal A}$:
$${\cal C}=\{(m,n,t)\in {\cal A}:|m|=2L\hbox{ or }|n|=2L\hbox{ or }t=0\}.$$

We define the model restricted to ${\cal A}+(kL,mL,nT)$ as the model with the same birth and death rates as the process on $\ZZ^2$ with the restriction that a pathogen in the complement of $[-2L,2L]^2+(kL,mL)$ cannot give birth inside $[-2L,2L]^2+(kL,mL)$ between times $nT$ and $(n+2)T$.

We will compare our models S2 and S3 to a percolation process on $\ZZ^2\times \ZZ_+$.
We declare $(k,m,n)$ in $\ZZ^2\times \ZZ_+$ to be wet if for the process restricted to ${\cal A}+(kL,mL,nT)$ there is no pathogen in ${\cal B}+(kL,mL,nT)$. Moreover, we want no  pathogen in ${\cal B}+(kL,mL,nT)$ for any possible configuration of the boundary ${\cal C}+(kL,mL,nT)$. 

Let $\epsilon>0$. We are going to show that given $\lambda>0$, there exists $r_1$ in (0,1) such that
$$P((k,m,n)\hbox{ is wet})\geq 1-\epsilon\hbox{ if }r<r_1.$$
Consider first models S2 and S3 restricted to ${\cal A}$ and with $r=0$. That is, a pathogen born inside ${\cal A}$ is always of the same type as its parent. Note that if there is a pathogen inside ${\cal B}$ there must be a line of infection from the boundary ${\cal C}$ of ${\cal A}$ to ${\cal B}$. This line of infection has either started at the bottom of ${\cal C}$ (i.e. at time 0) or on one its sides (i.e. at a time different from 0). We will now show that these possibilities have all exponentially small probability for both models S2 and S3. 

Since the line of infection cannot change type inside ${\cal A}$,
if there is a line of infection from the bottom of ${\cal C}$ to ${\cal B}$ then the type of the pathogens making up the line of infection must last at least $T$.  The death rate of a type is 1 for Model S2 and is at least 1 for Model S3. Hence, the probability that there is a line of infection from the bottom of ${\cal C}$ to ${\cal B}$ is less than $e^{-T}$. Note that there are $(4L+1)^2$ sites at the bottom of ${\cal C}$.

We now deal with a line of infection from a side of ${\cal C}$. The minimum distance between a side of ${\cal C}$ and ${\cal B}$ is $L$. Starting from a site $x$ on a side of ${\cal C}$ there are positive constants $c$, $C$ and $\gamma$ depending on $\lambda$ and such that the probability that a line of infection starting at $x$ reaches ${\cal B}$ by time $cL$ is less than
$Ce^{-\gamma L}$ (see for instance Lemma 9, p. 16 of Durrett (1988)). 
This estimate takes into account only births and so is valid for models S2 and S3.
If the line of infection takes at least $cL$ units of time to get to ${\cal B}$ then the type of the infection line from the side of ${\cal C}$ must last at least $cL$. For both models this has a probability less than $e^{-cL}$. Putting together these estimates we get
$$P((0,0,0)\hbox{ is wet})\geq 1-(4L+1)^2e^{-T}-8T(2L+1)Ce^{-\gamma L}-8T(2L+1)e^{-cL}\hbox{ for }r=0.$$
By taking $L=T$ large enough we get
$$P((0,0,0)\hbox{ is wet})\geq 1-\epsilon/2\hbox{ for }r=0.$$
Given that ${\cal A}$ is a finite box it is possible to find $r_1>0$ (depending on $\lambda$ and $\epsilon$) so close to 0  that if $r<r_1$ there is no creation of a new type in ${\cal A}$ with probability at least $1-\epsilon/2$. Thus, with probability at least $1-\epsilon/2$ models S2 and S3 restricted to ${\cal A}$ and with $r<r_1$ are coupled with models S2 and S3 with $r=0$, respectively. Hence,
$$P((0,0,0)\hbox{ is wet})\geq 1-\epsilon\hbox{ for }r<r_1.$$
By translation invariance, the same is true for any site $(k,m,n)$ of $\ZZ^2\times \ZZ_+$.
We now define a percolation process on $\ZZ^2\times \ZZ_+$ with finite range dependence. Let 
$${\cal A}(k,m,n)=(kL,mL,nT)+{\cal A}.$$
For each element $(k,m,n)$ in $\ZZ^2\times \ZZ_+$ we draw an oriented edge from $(k,m,n)$ to $(x,y,z)$ if $n\leq z$ and 
${\cal A}(k,m,n)\cap {\cal A}(x,y,z)\not=\emptyset.$
The wet sites in the ensuing directed graph constitute a percolation model.
The dependence of this percolation model has finite range because the event that $(k,m,n)$ is wet depends only on the Poisson processes inside ${\cal A}(k,m,n)$ and this box intersects only finitely many other boxes ${\cal A}(x,y,z)$. 
 
Note that if there is a pathogen somewhere at some time then there must be a path of dry sites in the percolation process. Since a dry (i.e. not wet) site has probability $\ee$ in this percolation process, 
by taking $\ee>0$ small enough one can make the probability of a path of dry sites of length $n$ decrease exponentially fast with $n$. This in turn implies that for any given site there will be no pathogen after a finite random time, see (8.2) in van den Berg et al. (1998).
\end{proof}

\noindent {\bf Remark.} This proof breaks down for Model S1 for at least two reasons. In the proof above it is crucial to have a lower bound on the death rate of a type. For Model S1 there is no such lower bound; the death rate of a type is $1/k$ if there are $k$ types and there is no upper bound on $k$. Moreover, what happens inside a finite space-time box for Model S1 depends on how many types there are in the whole space. Hence, there is little hope to compare Model S1 to a finite range percolation 
model as we did for models S2 and S3.

\begin{proof}[Proof of part 3 of Theorem \ref{spatial2thm}]
Let $e_1$ be the vector $(1,0,\dots,0)$ in $\ZZ^d$
$$B=[-2L,2L]^d\times[0,T]\qquad B_{m,n}=(4mLe_1,50nT)+B$$
$$I=[-J,J]^d$$
$${\cal L}=\{(m,n)\in \ZZ^2:m+n\hbox{ is even}\}.$$

We declare $(m,n)\in {\cal L}$ to be wet if there is $(x,t)$ in $B_{m,n}$ such that 
each site of the interval $x+I$ is occupied by a pathogen at time $t$
for the process restricted to $(4Lme_1,50nT)+(-6L,6L)^d\times [0,51T]$.

Set $r=1$ in the box $(-6L,6L)^d\times [0,51T]$.  
As noted before the set of occupied sites is a contact process for models S2 and S3. Since $\la>\la_c$ it is a supercritical contact process.

Bezuidenhout and Grimmett (1990) have shown that for a supercritical contact process, and for any $\epsilon>0$,
$J$, $L$ and $T$ can be chosen so that if $(0,0)$ is wet then with probability $1-\epsilon$, $(1,1)$ and
$(-1,1)$ will also be wet. Here we are following the approach and notation of Durrett (1991). 
More precisely, for any $\epsilon>0$ we can pick $J$, $L$ and $T$ such that 
$$P((1,1)\hbox{ and }(-1,1)\hbox{ are wet}|(0,0)\hbox{ is wet})>1-\epsilon\hbox{ for }r=1.$$
Since $(-6L,6L)^d\times [0,51T]$ is a finite space-time box, we can pick $r_2$ so close to 1 (but strictly smaller) that for models S2 and S3 with parameters $\la$ and $r>r_2$ all births inside $(-6L,6L)^d\times [0,51T]$ are of a new type with probability at least $1-\epsilon$. Therefore, the process of occupied sites for models S2 and S3 and $r>r_2$ may be coupled to a contact process with probability at least $1-\epsilon$.
Hence, for models S2 and S3 we have
$$P((1,1)\hbox{ and }(-1,1)\hbox{ are wet}|(0,0)\hbox{ is wet})>1-2\epsilon\hbox{ for }r>r_2.$$
By picking $\ee>0$ small enough we can show that there is a positive probability of an infinite wet cluster in ${\cal L}$. This, in turn, implies that pathogens have a positive probability of surviving forever, see Durrett (1991) for more details.
\end{proof}

\section{Analysis of the first spatial model}

In this section, we consider Model S1 and prove Theorem \ref{spatial1thm}, which says that the pathogens have positive probability of surviving whenever $\lambda > 0$ and $r > 0$.  This result is easiest to prove in $d \geq 2$, when there are always many neighboring sites on which the pathogens can give birth.  We begin by proving the result in this case.

\begin{proof}[Proof of Theorem \ref{spatial1thm} for $d \geq 2$]

Assume that at some time there are $n$ different types of pathogens in Model S1. Thus, there are $k\geq 
n$ occupied sites. It is easy to see that, if $d\geq 2$, at least ${\sqrt k}$ occupied sites have at least one empty neighbor. Therefore, the rate at which the number of types goes from $n$ to $n+1$ is at least $\lambda r\sqrt{k}\geq \lambda r\sqrt{n}$. On the other hand the rate at which the number of types goes from $n$ to $n-1$ is 1. Hence, the number of types in Model S1 is at least as large as a birth and death chain with transition rates: 
\begin{eqnarray*}
&& n\to n+1 \mbox{ at rate }\lambda r{\sqrt n}\\
&&n\to n-1 \mbox{ at rate } 1
\end{eqnarray*}

An argument very similar to the one in the proof of Theorem 1.1 shows that for all $\lambda > 0$ and $r > 0$, there is a positive probability that this chain never reaches zero. Since there are at least as many types in Model S1 as there are individuals in the birth and death chain, there is a positive probability that the number of types in S1 does not reach zero. This completes the proof of Theorem 2 for $d \geq 2$.
\end{proof}

\bigskip
We devote the rest of this section to the case $d = 1$.  This case is more complicated because if there are $n$ occupied sites, there could be as few as two sites with an empty neighbor.  Nevertheless, we will be able to show that the number of occupied sites grows linearly in time because the sites on the far left and the far right of the configuration will give birth at rate $\lambda$, while deaths create ``holes" in the configuration that usually fill up quickly.

We begin by introducing some notation.  Let $N(t)$ be the number of pathogens alive at time $t$, and let $N_k(t)$ be the number of type $k$ pathogens alive at time $t$.  Let $S_k$ be the set of sites that are occupied by a type $k$ pathogen at some time.  It is easy to see that $S_k$ is an interval.  Let $\zeta_k$ be the time at which the type $k$ pathogens die, with the convention that $\zeta_k = \infty$ if no type $k$ pathogen ever dies.  Define
\begin{align}
L(t) &= \inf\{x: \mbox{ there is a pathogen at site }x \mbox{ at time }t\}, \nonumber \\
R(t) &= \sup \{x: \mbox{ there is a pathogen at site }x \mbox{ at time }t\}, \nonumber
\end{align}
so all the pathogens at time $t$ are contained in the interval $[L(t), R(t)]$.  Fix a positive integer $T$.  Let $D(t)$ be the number of types that die before time $t$, and let $X(t)$ be the number of times $s$ with $T < s \leq t$ such that, at time $s$, either the type occupying the site $L(s-)$ or the type occupying the site $R(s-)$ dies at time $s$.  Let $Y(t)$ be the number of times $s \leq t$ such that, at time $s$, either the pathogen at site $L(s-)$ gives birth on site $L(s-) - 1$ or the pathogen on site $R(s-)$ gives birth on site $R(s-) + 1$.  Let $\zeta = \inf\{t: N(t) = 0\}$ be the time at which the pathogens die out, with the convention that $\zeta = \infty$ if $N(t) > 0$ for all $t$.  Let $\gamma = \min\{1, \lambda/6\}$, and let $\kappa = \inf\{t: N(t) < \gamma t\}$.

Assume that the initial configuration consists of $6T$ pathogens, all of different types labeled $1$, \dots, $6T$, with a pathogen at each of the sites $\{-3T + 1, \dots, 3T\}$.  Fix positive constants $C_1$, $C_2$, and $C_3$, and define the following six events:
\begin{itemize}
\item Let $A_1$ be the event that $D(t) \leq 2t$ for all $t > T$.

\item Let $A_2$ be the event that $Y(t) \leq 3 \lambda t$ for all $t > T$ and $Y(t) \geq Y(T) + \lambda(t - T)$ for all $2T < t < \zeta$.

\item Let $A_3$ be the event that for all $t > T$, we have ${\displaystyle \max_{0 \leq s \leq t} \max_k N_k(s)} \leq C_1 \log t$. 

\item Let $A_4$ be the event that for all $t > 2T$, at most $C_2 \log t$ different types die between times $t - C_3 \log t$ and $t$.

\item Let $A_5$ be the event that for all $2T < t < \kappa$, we have $X(t) \leq t^{1/2}$.

\item Let $A_6$ be the event that for all $2T < t < \kappa$, all $k \in {\bf N}$ such that $\zeta_k \leq t - C_3 \log t$, and all $x \in S_k$, there exists a time $s$ with $\zeta_k < s \leq t$ such that either $x < L(s)$, $x > R(s)$, or the site $x$ is occupied at time $s$.
\end{itemize}

\begin{Prop}
Let $\epsilon > 0$.  Then there exist positive constants $C_1$, $C_2$, and $C_3$ such that $$P \bigg( \bigcap_{i=1}^6 A_i \bigg) > 1 - 15 \epsilon$$ for sufficiently large $T$.
\label{6prop}
\end{Prop}

Before proving this proposition, we show that it implies Theorem \ref{spatial1thm} for $d = 1$.
The idea is that, on $A_2$, the length of the interval between the left-most pathogen and the right-most pathogen increases linearly.  On $A_3$, at most $C_1 \log t$ pathogens can die at time $t$, and on $A_1 \cap A_4$, deaths are sufficiently infrequent.  On $A_6$, when a type dies creating a ``hole" in the configuration, it fills up quickly, so that the interval between the left-most and right-most pathogens is mostly filled by pathogens.  This interval can get shorter when the left-most or right-most pathogen is killed, but on $A_5$ such deaths occur infrequently.

\begin{proof}[Proof of Theorem \ref{spatial1thm} for $d =1$]
If the process starts from a single pathogen at time zero, then with positive probability we eventually reach the configuration described above, with $6T$ pathogens of different types on sites $\{-3T + 1, \dots, 3T\}$.  Therefore, it suffices to show that if the process starts from a configuration with $6T$ pathogens of different types on $\{-3T + 1, \dots, 3T\}$, then with positive probability the process survives forever.  We will show that for sufficiently large $T$, we have $N(t) \geq \gamma t$ for all $t$ on the event $\cap_{i=1}^6 A_i$.  Theorem \ref{spatial1thm} for $d = 1$ will then follow from Proposition \ref{6prop}.

On $A_1$, we have $D(2T) \leq 4T$, so at least $2T$ of the $6T$ pathogens alive at time zero must survive until time $2T$.  Therefore, for all $t \leq 2T$, we have $N(t) \geq 2T \geq t \geq \gamma t$.  

Now assume $2T < t < \kappa$.  Suppose $L(t) < x < R(t)$ and there is no pathogen at $x$ at time $t$.  Because $L(t) < x < R(t)$, the site $x$ must have been occupied at some time before $t$.  Let $u = \sup\{s < t: \mbox{site }x \mbox{ is occupied at time }s\}$.  Suppose $u < t - C_3 \log t$.  Then on $A_6$, there is a time $s \in (u, t)$ such that either $x$ was occupied at time $s$, $x < L(s)$, or $x > R(s)$.  
However, since $L(t) < x < R(t)$, it follows that in the latter two cases, $x$ must be occupied at some time
in $(s, t)$, which contradicts the definition of $u$.  Thus, on $A_6$, we have $u \geq t - C_3 \log t$.
It follows that the number of vacant sites $x$ such that $L(t) < x < R(t)$ is at most the number of pathogens that die between times $t - C_3 \log t$ and $t$ which, on the event $A_3 \cap A_4$, is at most $C_1 C_2 (\log t)^2$.  This means that
\begin{equation}
N(t) \geq R(t) - L(t) - C_1 C_2 (\log t)^2.
\label{Nt}
\end{equation}

We now compare $Y(t)$ to $R(t) - L(t)$.  There are $Y(t) - Y(T)$ times at which the process $(R(s) - L(s), T \leq s \leq t)$ increases by one.  There are at most $X(t)$ times when the process decreases because
of deaths.  Suppose, at time $s \in (T, t]$, the pathogen at $L(s-)$ or $R(s-)$ dies.
Note that $(R(s-) - L(s-)) - (R(s) - L(s))$ is at most the number of pathogens that die at time $s$ plus the number of vacant sites between $L(s-)$ and $R(s-)$.  The number of pathogens that die is at most $C_1 \log t$ on $A_3$, and we just showed that the number of vacant sites between $L(s-)$ and $R(s-)$ is at most $C_1 C_2 (\log t)^2$ on $A_3 \cap A_4 \cap A_6$.  It follows that for $2T < t < \kappa$, we have
\begin{equation}
R(t) - L(t) \geq (R(T) - L(T)) + (Y(t) - Y(T)) - X(t)\big(C_1 \log t + C_1 C_2 (\log t)^2\big).
\label{RtLt}
\end{equation}
Now $X(t) \leq t^{1/2}$ on $A_5$.  On $A_2$, we have $Y(t) - Y(T) \geq \lambda (t - T) \geq \lambda t/2 \geq 3 \gamma t$.  Also, $R(T) - L(T) \geq 0$.  Therefore, by combining (\ref{Nt}) with (\ref{RtLt}), we see that on $\cap_{i=1}^6 A_i$,
$$N(t) \geq 3 \gamma t - t^{1/2}(C_1 \log t + C_1 C_2 (\log t)^2) - C_1 C_2 (\log t)^2.$$
It follows that for sufficiently large $T$, we have $N(t) \geq 2 \gamma t$ whenever $2T < t < \kappa$.

To show that $N(t) \geq \gamma t$ for all $t$ on $\cap_{i=1}^6 A_i$ for sufficiently large $T$, it remains to show that $\kappa = \infty$ on $\cap_{i=1}^6 A_i$ for sufficiently large $T$.  Because $N(t) \geq \gamma t$ for $t \leq 2T$ on $A_1$, we have $\kappa > 2T$.  Suppose $\kappa < \infty$.  On $A_3$, we
have $N(\kappa-) - N(\kappa) \leq C_1 \log \kappa$ because at most $C_1 \log \kappa$ pathogens can die at time $\kappa$.  However, on $\cap_{i=1}^6 A_i$, we have $N(\kappa-) \geq 2 \gamma \kappa$
and $N(\kappa) < \gamma \kappa$, so for $T$ large enough that $C_1 \log (2T) < 2 \gamma T$, we must have $\kappa = \infty$.
\end{proof}

Proposition \ref{6prop} will follow from Lemmas \ref{A1A2}, \ref{A3}, \ref{A4}, \ref{A5}, and \ref{A6} below.  Once $T$ is chosen sufficiently large, Lemma \ref{A1A2} will imply $P(A_1) > 1 - \epsilon$ and $P(A_2) > 1 - \epsilon$.  Lemma \ref{A3} then gives $P(A_3) > 1 - 3 \epsilon$, and Lemma \ref{A4} gives $P(A_4) > 1 - \epsilon$.  Finally, Lemma \ref{A5} implies $P(A_5) > 1 - 4 \epsilon$ and it follows from Lemma \ref{A6} that $P(A_6) > 1 - 5 \epsilon$.  Our first step will be to bound the probabilities of $A_1$ and $A_2$.

\begin{Lemma}
Let $\epsilon > 0$.  For sufficiently large $T$, we have $P(A_1) > 1 - \epsilon$ and $P(A_2) > 1 - \epsilon$.
\label{A1A2}
\end{Lemma}

\begin{proof}
Until time $\zeta$, deaths occur at times of a rate $1$ Poisson point process.  Let $D'(t)$ denote the number of points of a rate one Poisson process before time $t$, which can be coupled with the death process in such a way that $D(t) = D'(t)$ for all $t < \zeta$.  We have $t^{-1} D'(t) \rightarrow 1$ a.s.  It follows that $P(A_1) > 1 - \epsilon$ for sufficiently large $T$.

Likewise, until time $\zeta$, the pathogens at sites $L(t)$ and $R(t)$ each give birth on sites $L(t) - 1$ and $R(t) + 1$ respectively at rate $\lambda$, so these births occur at times of a Poisson point process of rate $2 \lambda$.  Let $Y'(t)$ denote the number of points of a rate $2 \lambda$ Poisson point process up to time $t$, coupled with the particle system so that $Y(t) = Y'(t)$ for $t < \zeta$.
Then $t^{-1} Y'(t) \rightarrow 2 \lambda$ a.s. and $(t - T)^{-1} (Y'(t) - Y'(T)) \rightarrow 2 \lambda$ a.s.  It follows that $P(A_2) > 1 - \epsilon$ for sufficiently large $T$.
\end{proof}

We next work towards bounding the probability of $A_3$.  The first step is to bound the probability that the number of pathogens of a given type is high.

\begin{Lemma}
Let $N_k(t)$ be the number of pathogens of type $k$ at time $t$.  Then there exist positive constants $C_4$ and $C_5$ such that for all $0 < r < 1$ and all $a$, we have $$P \big( \max_{t \geq 0} N_k(t) > a \big) \leq C_4 e^{-C_5 a}.$$
\label{typesize}
\end{Lemma}

\begin{proof}
It is clear from the description of the model that at any time $t$ at which there are pathogens of type $k$, the set of sites occupied by type $k$ pathogens is an interval of the form $\{a_t, a_t+1, \dots, b_t\}$.  The maximum number of type $k$ pathogens at any time can therefore be written as $1 + Y + Z$, where $Y$ is the number of times that the type $k$ pathogen on the far left of the interval gives birth on the site to its left, and $Z$ is the number of times that the type $k$ pathogen on the far right of the interval gives birth on the site to its right.  

Let $Z_1$ be the number of times that the type $k$ pathogen at site $b_t$ gives birth to another type $k$ pathogen on site $b_t + 1$, until the first time that the site $b_t+1$ is occupied by a pathogen of another type.  Because each pathogen born is a new type with probability $r$, the distribution of $Z_1 + 1$ is dominated by the geometric distribution with parameter $r$.  Once a different type occupies site $b_t+1$, the type $k$ pathogen at $b_t$ can not give birth again at site $b_t+1$ unless the type at $b_t+1$ dies before type $k$ dies, which happens with probability $1/2$.  It follows that the distribution of $Z$ is dominated by the distribution of $Z_1 + \dots + Z_N$, where $N$ has the geometric distribution with parameter $1/2$, $Z_i + 1$ has the geometric distribution with parameter $r$ for all $i$, and $N$ is independent of $Z_1, Z_2, \dots$.  Therefore,
\begin{align}
P \bigg( Z \geq \frac{a-1}{2} \bigg) &\leq \sum_{n=1}^{\infty} P(N = n) P \bigg( Z_i \geq \frac{a-1}{2n} \mbox{ for some }i \in \{1, \dots, n\} \bigg) \nonumber \\
&\leq \sum_{n=1}^{\infty} \frac{n}{2^n} (1 - r)^{((a-1)/2n) - 1} \leq (1-r)^{-3/2} \sum_{n=1}^{\infty} \frac{n}{2^{n}} (1-r)^{a/2n}.
\label{geombound}
\end{align}
In the sum on the right-hand side of (\ref{geombound}), the ratio of the $(n+1)$st term to the $n$th term converges to $1/2$ as $n \rightarrow \infty$, and therefore is less than $3/4$ for all $n \geq M$ for some integer $M$.  For this $M$, we have
$$P \bigg( Z \geq \frac{a-1}{2} \bigg) \leq (1-r)^{-3/2} \bigg( \frac{M}{2} (1-r)^{a/2M} + \frac{M}{2^M} \cdot \frac{(1-r)^{a/2M}}{1 - 3/4} \bigg) \leq \frac{C_4}{2} e^{-C_5 a},$$ where $C_4 = (1-r)^{-3/2} M (1 + 2^{3-M})$ and $C_5 = - \log(1-r)/2M$.  By the same argument, we get $P(Y \geq (a-1)/2) \leq (C_4/2) e^{-C_5 a}$.  Since the maximum number of type $k$ pathogens is $1 + Y + Z$, the result follows.
\end{proof}

\begin{Lemma}
Let $\epsilon > 0$.  There is a constant $C_1$ such that for sufficiently large $T$, we have $P(A_3^c \cap A_1 \cap A_2) < \epsilon$.
\label{A3}
\end{Lemma}

\begin{proof}
Suppose $t > T$.  On $A_2$, there is a set of at most $3 \lambda t + 6T \leq 3 (\lambda + 2) t$ sites 
at which there has been a pathogen at some time $s \leq t$.  Since there are at most $2t$ deaths before time $t$ on $A_1$, no site can be occupied by more than $2t + 1$ different pathogens before time $t$. 
Therefore, on $A_1 \cap A_2$, at most $3(\lambda + 2)(2t+1)t$ different types of pathogens can be born by time $t$.

For positive integers $n$, let $t_n = T^{2^n}$.  Given a constant $C_1$, let $B_n$ be the event that for some $k \leq 3(\lambda + 2)(2t_n + 1)t_n$, we have $N_k(s) > \frac{1}{2} C_1 \log t_n$ for some $s$. 
On the event $A_3^c \cap A_1 \cap A_2$, there is a $t > T$ and a $k \leq 3 (\lambda + 2)(2t + 1)t$ such that $N_k(s) > C_1 \log t$ for some $s$.  If $n = \inf\{m: t_m \geq t\}$, then $k \leq 3 (\lambda + 2)(2t_n + 1)t_n$ and $N_k(s) > C_1 \log t > C_1 \log t_{n-1} = \frac{1}{2} C_1 \log t_n$.
Therefore, on $A_3^c \cap A_1 \cap A_2$, the event $B_n$ occurs for some $n$.  By Lemma \ref{typesize}, we have $$P(B_n) \leq 3(\lambda + 2)(2t_n + 1)t_n \cdot C_4 e^{-C_1 C_5 (\log t_n)/2} = 3C_4(\lambda + 2)(2t_n + 1)t_n^{1 - C_1C_5/2}.$$
Therefore, $$P(A_3^c \cap A_1 \cap A_2) \leq 3C_4(\lambda + 2) \sum_{n=1}^{\infty} (2t_n + 1)t_n^{1-C_1C_5/2}.$$  If we choose $C_1$ large enough that $2 - C_1C_5/2 < 0$, then this expression is less than $\epsilon$ for sufficiently large $T$.
\end{proof}

The next two results bound the probabilities of $A_4$ and $A_5$.  Both of these events pertain to the number of deaths.  The proofs make use of the fact that, up to time $\zeta$, deaths occur at times of a rate one Poisson process.
We first state a lemma related to the gamma distribution, which we can use to choose the constant $C_3$, now that $C_1$ has already been chosen.  The reason for this choice will become clear later.  We will then choose $C_2$ in Lemma \ref{A4}.

\begin{Lemma}
There exists a constant $C_3$ such that if $X$ has a gamma distribution with shape parameter $2 C_1 \log (t + C_3 \log t)$ and scale parameter $2 \lambda$, then $P(X > C_3 \log t) \leq t^{-2}$ for sufficiently large $t$.
\label{gamlem}
\end{Lemma}

\begin{proof}
If $0 < \theta < 2 \lambda$, then $E[e^{\theta X}] = [2 \lambda/(2 \lambda - \theta)]^{2 C_1 \log (t + C_3 \log t)}$.  It follows from Markov's Inequality, taking $\theta = \lambda$, that $$P(X > C_3 \log t) \leq e^{-\lambda C_3 \log t} E[e^{\lambda X}] = e^{-\lambda C_3 \log t} 2^{2 C_1 \log (t + C_3 \log t)}.$$
Let $g(t) = \log(t + C_3 \log t)/(\log t)$.  Then $$P(X > C_3 \log t) \leq t^{-\lambda C_3 + 2 C_1 g(t)}.$$  We can choose $C_3$ such that $\lambda C_3 - 4C_1 \geq 2$ and then $t$ large enough that $g(t) \leq 2$.  The lemma follows.
\end{proof}

\begin{Lemma}
Let $\epsilon > 0$.  There is a constant $C_2$ such that $P(A_4) > 1 - \epsilon$ for sufficiently large $T$.
\label{A4}
\end{Lemma}

\begin{proof}
For integers $n \geq 1$ and $k \geq 0$ such that $2^nT + C_3 (k-1) \log (2^{n+1}T) < 2^{n+1}T$,
define the interval $$I_{n,k} = \big[ 2^nT + C_3 (k-1) \log(2^{n+1}T), \: 2^nT + C_3 k \log(2^{n+1}T) \big].$$
Let $B_{n,k}$ be the event that at least $\frac{1}{2} C_2 \log (2^n T)$ types die during the time interval 
$I_{n,k}$.  For any $t$ such that $2^nT \leq t < 2^{n+1}T$, the interval from $t - C_3 \log t$ to $t$ is contained in $I_{n,k} \cup I_{n, k+1}$ for some $k$.  Therefore, if for some $t$ such that $2^nT \leq t < 2^{n+1}T$, more than $C_2 \log t$ types die between times $t - C_3 \log t$ and $t$, the event $B_{n,k}$ must occur for some $k$.  It follows that $$P(A_4^c) \leq \sum_{n=1}^{\infty} \sum_k P(B_{n,k}),$$ so we need to bound the probabilities $P(B_{n,k})$.

Until time $\zeta$, types die at times of a rate one Poisson process, so the distribution of the number of types that die during the interval $I_{n,k}$ is dominated by the Poisson distribution with mean $C_3 \log (2^{n+1}T)$.  If $X$ has the Poisson distribution with mean $\lambda$, then for all $\theta > 0$, we have $$P(X \geq a \lambda) \leq e^{-\theta a \lambda } E[e^{\theta X}] = e^{-\theta a \lambda + \lambda(e^{\theta} - 1)}.$$  Choosing $\theta = \log a$, we get
\begin{equation}
P(X \geq a \lambda) \leq e^{-\lambda (a \log a - a + 1)}.
\label{Poisbound}
\end{equation}
To bound $P(B_{n,k})$, we need to apply $(\ref{Poisbound})$ with $\lambda = C_3 \log (2^{n+1} T)$ and $a \lambda = \frac{1}{2} C_2 \log (2^n T)$.  This means that $a = (C_2 \log (2^n T))/(2 C_3 \log (2^{n+1} T))$.  We can choose $C_2$ large enough that, for all $n$, we have $b = C_3 (a \log a - a + 1) > 1$.  For this choice of $C_2$, we get $$P(B_{n,k}) \leq e^{-b \log(2^{n+1} T)} = (2^{n+1}T)^{-b}.$$
For sufficiently large $T$, we have that for all $n$ there are at most $2^{n+1}T$ intervals $I_{n,k}$.  For such $T$, $$P(A_4^c) \leq \sum_{n=1}^{\infty} (2^{n+1}T)^{1-b},$$ which is less than $\epsilon$ for sufficiently large $T$.
\end{proof}

\begin{Lemma}
Let $\epsilon > 0$.  For sufficiently large $T$, we have $P(A_5^c \cap A_3) < \epsilon$.
\label{A5}
\end{Lemma}

\begin{proof}
Until time $\zeta$, deaths occur at times of a rate one Poisson process.  Denote the times of such a Poisson process by $0 < \tau_1 < \tau_2 < \dots$.  Define a sequence of independent random variables $(U_i)_{i=1}^{\infty}$, each having a uniform distribution on $[0, 1]$.  When a death event occurs, one type is chosen at random to die.  Therefore, denoting the number of types at time $t$ by $M(t)$, we may assume that until time $\zeta$, deaths occur at the times $\tau_1 < \tau_2 < \dots$ and that, at time $\tau_i$, if $M(\tau_i-) \geq 2$ then either the type at $L(\tau_i-)$ or $R(\tau_i-)$ dies if and only if $U_i \leq 2/M(t)$.

Suppose $T < t < \kappa$.  On the event $A_3$, the number of pathogens of a given type before time $t$ is at most $C_1 \log t$, so the number of types is at least $N(t)/(C_1 \log t) \geq \gamma t/(C_1 \log t)$.  
It follows that $X(t)$ is at most the number of times $\tau_i$ such that $T < \tau_i \leq t$ and $U_i \leq (2 C_1 \log \tau_i)/(\gamma \tau_i)$.  Such times $\tau_i$ occur at times of an inhomogeneous Poisson process of rate $\lambda(s) = (2 C_1 \log s)/(\gamma s)$.  It follows that $P(A_5^c \cap A_3)$ is at most the probability that, for some $t$, such a Poisson process contains at least $t^{1/2}$ points between times $T$ and $t$.  

For positive integers $n$, let $t_n = 4^n T$.  Let $B_n$ be the event that there are at least $\frac{1}{2} t_n^{1/2}$ points of the Poisson process between times $T$ and $t_n$.  If there is a $t$ such that there are at least $t^{1/2}$ points between times $T$ and $t$, then if $n = \min\{m: 4^m T \geq t\}$, there are at least $t^{1/2} \geq (t_n/4)^{1/2} = \frac{1}{2}t_n^{1/2}$ points between times $T$ and $t_n$, so $B_n$ occurs.  It follows that $P(A_5^c \cap A_3) \leq \sum_{n=1}^{\infty} P(B_n)$.

To bound $P(B_n)$, first note that for $T \geq 1$, the distribution of the number of points of the Poisson process between times $T$ and $t_n$ is Poisson with mean
$$\int_T^{t_n} \frac{2 C_1 \log s}{\gamma s} \: ds \leq \frac{2 C_1 \log t_n}{\gamma} \int_1^{t_n} \frac{1}{s} \: ds = \frac{2 C_1 (\log t_n)^2}{\gamma}.$$
Thus, $P(B_n)$ is at most the right-hand side of (\ref{Poisbound}) when $\lambda = 2 \gamma^{-1} C_1 (\log t_n)^2$ and $a \lambda = \frac{1}{2} t_n^{1/2}$.  It follows easily that $\sum_{n=1}^{\infty} P(B_n) < \epsilon$ for sufficiently large $T$, which completes the proof.
\end{proof}

It remains to bound $P(A_6)$.  Informally, $A_6$ is the event that whenever a type dies, creating a ``hole" in the configuration, the hole fills up quickly.  

Suppose a type dies at time $t$.  The pathogens that died occupied some interval $[\ell_t, r_t]$.  If the sites $\ell_t - 1$ and $r_t + 1$ are occupied at time $t$, then we say that a hole is created at time $t$.  
If $s > t$, we say that the hole exists until time $s$ if there is a $\ZZ$-valued process $(H(u), t \leq u \leq s)$ such that:
\begin{itemize}
\item $H(t) \in [\ell_t, r_t]$.

\item The site $H(u)$ is empty for all $u \in [t,s]$.

\item $H(u-) - 1 \leq H(u) \leq H(u-) + 1$ for all $u \in [t, s]$.

\item $L(u) < H(u) < R(u)$ for all $u \in [t, s]$.
\end{itemize}
For each $u \in [t, s]$, we think of $H(u)$ as being a site in the hole that was created at time $t$.  Over time, this site may move around within the hole so that no pathogen is born on it.  Note that it need not be the case that $H(u) \in [\ell_t, r_t]$ for all $u \in [t, s]$.  For example, if the pathogen occupying site $r_t + 1$ dies, the hole can exist beyond the time at which pathogens are born on all of the sites in $[\ell_t, r_t]$ if the site $r_t + 1$ remains vacant.
If no such process $(H(u), t \leq u \leq s)$ exists, then we say the hole disappears by time $s$.

\begin{Lemma}
Let $\epsilon > 0$.  For sufficiently large $T$, we have $P(A_6^c \cap A_1 \cap A_3) < \epsilon$.
\label{A6}
\end{Lemma}

\begin{proof}
We call a hole long-lasting if it is created at time $t \leq 2T - C_3 \log (2T)$ and exists until time $t + C_3 \log (2T)$, or if it is created at time $t \geq 2T - C_3 \log (2T)$ and exists until time $t + C_3 \log t$.  Assume $T$ is large enough that $2T - C_3 \log (2T) \geq T$ and that the function $t \mapsto t - C_3 \log t$ is increasing on $[T, \infty)$.  If $2T < t < \kappa$ and there are no long-lasting holes created before 
time $\kappa - C_3 \log \kappa$, then every hole created before time $t - C_3 \log t$ disappears by time $$\max\{2T, t - C_3 \log t + C_3 \log (t - C_3 \log t)\} \leq t.$$  
It follows that if there are no long-lasting holes created before time $\kappa - C_3 \log \kappa$, then $A_6$ must occur.  This is because if $A_6$ does not occur, then there exist $t \in (2T, \kappa)$, $\zeta_k \leq t - C_3 \log t$, and $x \in S_k$ such that for all $s \in [\zeta_k, t]$, the site $x$ is vacant at time $s$ and $L(s) < x < R(s)$.  Therefore, if a new hole is created at time $\zeta_k$, then by taking $H(s) = x$ for all $s \in [\zeta_k, t]$, we see that the hole exists until time $t$, contradicting that there are no long-lasting holes.  If a new hole was not created at time $\zeta_k$, then a hole created at an even earlier time lasts until time $t$, which gives the same contradiction.  It thus remains to bound the probability that there is such a long-lasting hole.

Suppose a hole is created at time $t$.  The pathogens that died at time $t$ occupied some interval $[\ell_t, r_t]$.  Label ``$a$" the type occupying $\ell_t - 1$ at time $t$, and label ``$b$" the type occupying $r_t + 1$ at time $t$.  Label ``$c$" the type at site $\max\{x < \ell_t:$ there is not a pathogen of type $a$ at site $x\}$, if this site is occupied.  Likewise, label ``$d$" the type at the site $\min\{x > r_t:$ there is not a pathogen of type $b$ at site $x\}$, if this site is occupied.  

Suppose $t \leq 2T - C_3 \log (2T)$, and that no holes created at earlier times are long-lasting. 
As long as a hole exists, pathogens are giving birth at rate $\lambda$ on the sites on the endpoints of the hole.  Therefore, the size of the hole decreases by one at times of a rate $2 \lambda$ Poisson process, until the hole no longer exists.  By Lemma \ref{gamlem}, the probability that fewer than $2 C_1 \log (2T)$ points of this Poisson process occur by time $C_3 \log (2T)$ is at most $(2T)^{-2}$ for sufficiently large $T$.  On $A_3$, before time $2T$ there can be no more than $C_1 \log (2T)$ pathogens of a given type.  Therefore, if more than $2 C_1 \log (2T)$ points of the Poisson process occur by time $C_3 \log (2T)$, the hole will not exist at time $t + C_3 \log (2T)$ unless two of the types $a$, $b$, $c$, and $d$ die before time $t + C_3 \log (2T)$.  This is clear if all four types exist.  If, for example, type $c$ does not exist, then if type $a$ dies at time $t^*$, either the hole will not exist beyond time $t^*$ because $L(t^*)$ will be to the right of where the hole was previously, or the hole will merge with a hole born before time $t$, which by assumption is not long-lasting and therefore will not still exist at time $t + C_3 \log (2T)$.  However, before time $2T$, there are always at least $2T$ types on $A_1$, so the rate at which one of these four types is dying is at most $4/(2T)$.  Using that when $X$ has a Poisson distribution with parameter $\lambda$, we have $P(X \geq 2) \leq \lambda^2$, we see that the probability that two of the four types die by time $C_3 \log (2T)$ is at most $[4 C_3 (\log 2T)/2T]^2$, and therefore the probability that the hole created at time $t$ is long-lasting is at most $[(4 C_3 (\log 2T) + 1)/2T]^2$.

Suppose instead $2T - C_3 \log (2T) \leq t \leq \kappa - C_3 \log \kappa$, and that no holes created at earlier times are long-lasting.  On $A_3$, before time $t + C_3 \log t$ there can be no more than $C_1 \log(t + C_3 \log t)$ pathogens of a given type.  Therefore, the hole can be long-lasting only if either two of the types $a$, $b$, $c$, and $d$ die before time $t + C_3 \log t$, or if there are fewer than $2 C_1 \log(t + C_3 \log t)$ points of a rate $2 \lambda$ Poisson process (whose points correspond to births at the endpoints of the hole) before time $C_3 \log t$.  Lemma \ref{gamlem} implies that the probability of the latter is at most $t^{-2}$ for sufficiently large $T$.  Since $t + C_3 \log t < \kappa$, the number of types during the interval from $t$ to $t + C_3 \log t$ is always at least $\gamma t/(C_1 \log(t + C_3 \log t))$.  Therefore, the rate of deaths of the four types is at most $4 C_1 \log(t + C_3 \log t)/\gamma t$, so the probability of at least two deaths during this time interval is at most $[4 C_1 C_3 (\log t) \log(t + C_3 \log t)/\gamma t]^2$.  Therefore, the probability that the hole created at time $t$ is long-lasting is at most $[(4 C_1 C_3 (\log t) \log(t + C_3 \log t) + \gamma)/\gamma t]^2$.

Since deaths occur at rate $1$, the bounds in these two time intervals imply that
\begin{align}
P(A_6^c \cap A_1 \cap A_3) &\leq \int_0^{2T - C_3 \log (2T)} \bigg( \frac{4 C_3 (\log 2T) + 1}{2T} \bigg)^2 \: dt \nonumber \\
&\hspace{.3in} + \int_{2T - C_3 \log (2T)}^{\infty} \bigg( \frac{4 C_1 C_3 (\log t) \log(t + C_3 \log t) + \gamma}{\gamma t} \bigg)^2 \: dt, \nonumber
\end{align}
which is less than $\epsilon$ for sufficiently large $T$.
\end{proof}

\bigskip
\begin{center}
{\bf {\Large Acknowledgments}}
\end{center}

\noindent We thank two referees for helpful comments and corrections.

\bigskip
\begin{center}
{\bf {\Large References}}
\end{center}

\mn M. Bernaschi and F. Castiglione (2002). Selection of escape mutants from immune recognition during HIV infection. {\it Immunology and Cell Biology}, {\bf 80}, 307-313.

\mn C. Bezuidenhout, G. Grimmett (1990). The critical contact process dies out.  {\it Ann. Probab.} {\bf 18}, 1462-1482.

\mn F. Chung and L. Lu (2004).  Coupling online and offline analyses for random power law graphs.  {\it Internet Mathematics} {\bf 1}, 409-461.

\mn C. Cooper, A. Frieze, and J. Vera (2004).  Random deletion in a scale-free random graph process.  {\it Internet Mathematics} {\bf 1}, 463-483.

\mn R. J. De Boer and A. S. Perelson (1998). Target cell limited and immune control models of HIV infection: a comparison. {\it Journal of Theoretical Biology}, {\bf 190}, 201-214.

\mn R. Durrett (1988). {\it Lecture notes on particle systems and percolation.} Wadsworth, Pacific Grove, California.

\mn R. Durrett (1991). The contact process, 1974-1989. {\it Lectures in Applied Mathematics}, 27, 1-18,
American Mathematical society.

\mn R. Durrett and J. Schweinsberg (2005).  Power laws for family sizes in a duplication model.  To appear in {\it Ann. Probab.}  Preprint available at http://front.math.ucdavis.edu/math.PR/0406216.

\mn P. G. Hoel, S. C. Stone and C. J. Stone (1972). {\it Introduction to Stochastic Processes.} Houghton Mifflin Company Boston.

\mn Y. Iwasa, F. Michor and M. Nowak (2004). Some basic properties of immune selection. {\it Journal of Theoretical Biology}, {\bf 229}, 179-188.

\mn T. Liggett (1999). {\it Stochastic Interacting Systems: Contact, Voter and Exclusion Processes}, Springer, Berlin.

\mn M. A. Nowak and R. M. May (2000). {\it Virus Dynamics, Mathematical Principles of Immunology and Virology.} Oxford University Press, Oxford.

\mn A. S. Perelson and G. Weisbuch (1997). Immunology for physicists. {\it Rev. Mod. Phys.}, {\bf 69}, 1219-1267.

\mn A. Sasaki (1994). Evolution of Antigen Drift/Switching: continuously evading pathogens. {\it Journal of Theoretical Biology}, {\bf 168}, 291-308. 

\mn R. B. Schinazi (1999). {\it Classical and spatial stochastic processes.} Birkhauser, Boston.

\mn G. Silvestri and M.B. Feinberg (2003). Turnover of lymphocytes and conceptual paradigms in HIV infection. {The Journal of Clinical Investigation}, {\bf 112}, 821-823.

\mn J. van den Berg, G. Grimmett, and R. Schinazi (1998). Dependent random graphs 
and spatial epidemics.  {\it Ann. Appl. Probab.} {\bf 8}, 317-336.

\mn G. U. Yule (1925).  A mathematical theory of evolution based on the
conclusions of Dr. J. C. Willis. {\it Phil. Trans. Roy. Soc. London Series B}.
{\bf 213}, 21--87.

\mn R. M. Zorzenon dos Santos (1999). Immune responses: getting close to experimental results with cellular automata models, Stauffer D. editor, {\it Ann. Rev. Comp. Phys. VI}, Singapore, World Scientific. 

\end{document}